\documentclass[]{svjour3}
\usepackage{amssymb} 
\usepackage{bm}
\usepackage{amsmath}
\usepackage{graphicx}

\newcommand\F{\mbox{I\kern-2pt F}}

\newcommand\cF{{\mathcal F}}

\newcommand\cH{{\mathcal H}}

\newcommand\cB{{\mathcal B}}

\newcommand\cU{{\mathcal U}}

\newcommand\R{{\mathbb R}}

\def\text#1{\hbox{#1}}

\def\E{{\bf E}}
\def\P{{\bf P}}

\def\build #1_#2{\mathrel{\mathop{\kern 0pt #1}\limits_{#2}}} 

\def\R{{\mathbb R}}

\newcommand\fdem{$\Box$}

\newcommand\beq{\begin{equation}}
\newcommand\eeq{\end{equation}}
\newcommand\bea{\begin{eqnarray}}
\newcommand\eea{\end{eqnarray}}
\newcommand\bean{\begin{eqnarray*}}
\newcommand\eean{\end{eqnarray*}}

\newtheorem{ssmptn}[theorem]{Assumption}

\begin{document}

\title{Seierstad Sufficient Conditions for Stochastic Optimal Control Problems with Infinite Horizon}


\author{Anton\,O.\,Belyakov, Yuri\,M.\,Kabanov, Ivan\,A.\,Terekhov, Maxim\,M.\,Savinov}

\institute{Anton\,O.\,Belyakov,  Corresponding author  \at
              Lomonosov Moscow State University\\
              National University of Science and Technology ``MISiS''\\
              Moscow, Russia\\
              belyakov@mse-msu.ru
           \and
              Yuri\,M.\, Kabanov \at
              LMB, Universit\'e de Franche-Comt\'e \\
              Besan\c{c}on, France\\
              Lomonosov Moscow State University\\
              Moscow, Russia
	\and
	   Ivan\,A.\,Terekhov, Maxim\,M.\,Savinov \at
              Lomonosov Moscow State University\\
              ``Vega" Institute\\
              Moscow, Russia
}

\date{}

\maketitle

\begin{abstract}
In this note we consider a problem of stochastic  optimal control with the infinite-time horizon. We 
present analogues of the Seierstad sufficient conditions of overtaking optimality based on 
the dual variables stochastic described by BSDEs appeared 
in the Bismut--Pontryagin maximum principle.  
\end{abstract}


\section{Introduction}

Infinite-horizon optimal control problems play an important role in the modern economic theory.  A typical examples  are  optimal investment problems and optimal consumptions of households. In both cases  stochastic parameters are natural, see, e.g., \cite{Romer}. Historically, the main tools for  stochastic optimal control problems were HJB equations --- in a striking contrast 
with the deterministic optimal control problems, where many results were obtained using the Pontryagin maximum principle. 

Necessary optimality conditions in the form of the stochastic maximum principle are discussed already for several decades. They appeared via convex duality in the seminal paper by Bismut \cite{Bismut1973}, where   backward stochastic differential equations (BSDEs) with linear driver  were introduced  to describe the dynamics of dual variables. Bismut applied successfully his technique to the stochastic linear-quadratic regulator and the Merton consumption-investments 
problems in \cite{Bismut1976a,Bismut1976}. By this reason we believe that his contribution merits to be called the Bismut--Pontryagin stochastic maximum principle.  The relation between BSDEs  and
 the predictable representation theorem for the Wiener functionals was understood  
in the rarely available paper   \cite{Kabanov1978}, where stochastic maximum principle in the form of necessary and sufficient condition was established for a model with linear dynamics of the phase variable. Arkin and Saksonov  discussed in \cite{ArkinSaksonov1978}  some  important features of the stochastic maximum principle in the non-convex case.  Nowadays the mentioned papers have only a historical interest. The modern theory of non-linear BSDEs was developed  in much deeper works   by Pardoux and Peng, \cite{PP1990,P1990}. The literature in this field and its applications  is still exploding,  see, e.g., a more recent book \cite{CZ} on stochastic contract theory or \cite{OksendalSulem} on stochastic control of jump diffusions.      

%

In this paper we suggest  a stochastic analogue of the Seierstad theorem, \cite{Seierstad1977,Seierstad1978,Belyakov2020}, on sufficient 
conditions for overtaking optimality.  As in the Mangasarian conditions, it is assumed that the Hamilton--Pontryagin function is concave in the state and control variables. The obtained sufficient conditions are written using solutions of the BSDEs depending on the running time horizon.  

We also prove that in the case of a linear Hamilton--Pontryagin function the Seierstad conditions become not only sufficient but also necessary for the overtaking optimality of the control.

We provide examples, where the optimality conditions can be checked. In particular, in Example 2 we consider a stochastic version of the Nerlove--Arrow dynamics, see \cite{NA,Marinelli} and 
present a result on the overtaking optimality for the optimal investment problem  with quadratic investment adaptation costs and stochastic marginal productivity of capital, cf.  \cite{Romer}.

\section{Setting and main result}
\label{s:sp}
Let $(\Omega,\cF,{\bf F}=(\cF_t)_{t\ge 0},\P)$ be a filtered probability space satisfy the usual condition. Assume that ${\bf F}$  is generated by a Wiener process ${\bf  W}=(W_t)_{t\ge 0}$.  

%

Let $L^{22}_T:=L^2(\Omega\times [0,T],\cF\otimes \cB[0,T], dPdt)$. 

Fix a non-empty subset $U$ of  $\R$.  The control process ${\bf u}=(u_t)_{t\ge 0}$ in our setting is just a $U$-valued predictable process. We say that  the control ${\bf u}$ is admissible if  $(u_t)_{t\in [0,T]}\in L^{22}_T$ for each $T$. We denote $\cU$ the set of admissible controls.

The controlled
  dynamics corresponding to ${\bf u}$ is  given by the scalar process  $\boldsymbol x=\boldsymbol x^{\boldsymbol u}=(x_t)$ with a fixed initial value $x_0$ solving the  
stochastic differential equation 
$$
    dx_t= f(t, u_t,x_t)\,dt+\sigma (t, u_t,x_t)\, dW_t. 
$$

To guarantee the existence and uniqueness of strong solution we  assume  several rather standard hypotheses. 

The drift coefficient  $f=f(\omega,t,u,x)$ is continuously differentiable  in $(u,x)$ for each $(\omega,t)$ and predictable in $(\omega,t)$ for each 
 $(u,x)\in U\times \R$, that is the function  $(\omega,t)\mapsto f(\omega,t,u,x)$ is predictable.  The derivative in $x$  is uniformly bounded:  $|f_x(\omega,t,u,x)|\le {\rm const}$.  Moreover,  the  
  linear growth condition is fulfilled: 
$$
|f(\omega,t,u,x)|\le C_t(1+|u|+|x|), \qquad \forall\, (\omega,t,u)\in \Omega\times \R_+\times U\times \R, 
$$
where  ${\bf C}=(C_t)$ is a finite increasing  function

The  diffusion coefficient $\sigma =\sigma(\omega,t,u,x)$ satisfy the same conditions as the drift coefficient.   
It follows  that 
$$
\E \int_0^T \sigma^2(t,u_t,x_t)dt<\infty. 
$$

For a given finite time horizon $T$,   we consider   the value functional (to maximize)
\begin{equation}\label{eq:J}
  J_T(\boldsymbol{u}) := \E \int_0^T g(t, u_t, x_t) dt,
\end{equation}  
where the process $x={ \bf x}^{ \bf u}$.

We assume that the function  $g=g(\omega,t,u,x)$ is  continuously differentiable  in $(u,x)$ for each $(\omega,t)$ and predictable in $(\omega,t)$ for each 
 $(u,x)\in U\times \R$ and, for simplicity, satisfies the linear growth condition as $f$ and $\sigma$. 
 Moreover, we assume that 
$$
\big|g_x(\omega,t,u,x)\big|\le  \tilde c_T(1+|u|+|x|), \qquad  \forall\,(\omega,t,u)\in \Omega\times [0,t]\times U\times \R,
$$ 
where  $  \tilde{\bf c}=(\tilde c_t)$ is a finite increasing  function. This condition ensures that for any admissible control 
\smallskip
\begin{equation}
\label{gx}
    {\bf E} \int_0^T  \left| g_x(t, u_t, x_t)\right|^2 dt<\infty.
\end{equation}
These assumptions imply that for any admissible control    $J_T(\boldsymbol{u})<\infty$.
 
 However, the linear growth assumption of $g$ is too restrictive (it excludes  quadratic cost). 
 One can work without restricting  controls to those with well-defined value functional.

\begin{definition} The  control $\hat{\bf u}\in \cU$ is called {\it weakly overtaking optimal, WOO}, if 
for any   $ {\bf u}\in \cU$
$$
\liminf_{T \to \infty} \left(J_T(\boldsymbol{u})-J_T(\hat{\boldsymbol{u}}) \right) \le 0.
$$ 
\end{definition}

\begin{definition}  The  control $\hat{\bf u}\in \cU$ is called {\it overtaking optimal, OO}, if 
for any   $ {\bf u}\in \cU$
$$
\limsup_{T \to \infty} \left(J_T(\boldsymbol{u})-J_T(\hat{\boldsymbol{u}}) \right) \le 0. 
$$ 
\end{definition}

Let us define the Hamilton--Pontryagin function
 \begin{equation}
 \label{eq5}
 \mathcal{H}(t,u, x, p, h)=p\,f(t,u,x)+h\,\sigma(t,u,x) + g(t, u, x).
\end{equation}
We skip, as usual,   $\omega$ in this formula. 

Let $\hat{\bf u}$ be an admissible control which we want to test for optimality and let $\hat {\bf x}$ be the corresponding process.  
We consider on the interval  $[0,T]$ the BSDE
\begin{align}\label{eq6}
    dp_t=-\mathcal{H}_x(t, \hat u_t, \hat x_t, p_t,h_t)dt + h_tdW_t,  \quad p_T=0.
\end{align}
The driver of this  equation, 
$$
\mathcal{H}_x(t, \hat u_t(\omega), \hat x_t(\omega), p,h):=pf_x(t,\hat  u_t, \hat x_t)+h\sigma _x(t,\hat  u_t, \hat x_t)+g_x(t,\hat  u_t, \hat x_t),
$$
is such that  
$$
\E \int_0^T |\mathcal{H}_x(t, \hat u_t, \hat x_t, 0,0)|^2dt =\E \int_0^T|g_x(t,\hat  u_t, \hat x_t)|^2dt<\infty 
$$
in virtue of (\ref{gx}). Also, 
\begin{align*}
&|\mathcal{H}_x(t, \hat u_t(\omega), \hat x_t(\omega), p,h)-\mathcal{H}_x(t, \hat u_t(\omega), \hat x_t(\omega), \tilde p,\tilde h)|\\&\le
|p-\tilde p|f_x(t,\hat  u_t, \hat x_t)+|h-\tilde h|\sigma_x(t,\hat  u_t, \hat x_t)\\
&\le C(|p-\tilde p|+|h-\tilde h|).
\end{align*}
These conditions ensure
the existence and uniqueness of solution 
$({\bf p}^T,{\bf h}^T)$ of BSDE (\ref{eq6}) in the space $X_T\times L^{22}_T$, where $X_T$ is the space of adapted continuous processes ${\bf p}$ such that $\E\sup_{s\le T}|p_s|^2<\infty$, see \cite{CZ}.   
  

%



Our hypotheses  guarantee the existence and uniqueness of  solution   $(\hat{\bf p}, \hat{\bf h})=(\hat{\bf p}^T, \hat{\bf h}^T)$
in the Banach space $B_T\times L^{22}_T$, where $B_T$ is the space of continuous adapted processes  $p$ such that $\E \sup_{t\le T}
|p_t|^2<\infty$.


\begin{ssmptn}\label{a:concave}
There is  a convex set $\tilde U\supseteq U$ and  $T_1>0$ such that for  all $t \ge T_1$ the  function
$(u,x)\mapsto \mathcal{H}(t, u, x,p,h)$ is concave on  $ \tilde{U} \times \R$.
\end{ssmptn} 

\begin{theorem}(Optimality: sufficient conditions.)
\label{th:S}
  Let Assumption  \ref{a:concave}  be fulfilled. Then the admissible control  $(\hat{\boldsymbol{u}})$ is WOO, 
  if
    \begin{equation}
    \label{eq:dHsup}
    \liminf_{T \to \infty} \E \int_{0}^{T} \cH_u(t,\hat u_t, \hat x_t, p_t^T,h^T_t) \left( u_t - \hat u_t\right)  dt \le 0,
\end{equation}
 and OO,   if
\begin{equation}
\label{eq:dHinf}
       \limsup_{T \to \infty} \E \int_{0}^{T} \cH_u(t, \hat u_t, \hat x_t, p_t^T,h^T_t) \left(u_t -\hat  u_t\right)  dt \le 0,
  \end{equation}
 for every admissible control $\boldsymbol{u}$.
\end{theorem}

\begin{theorem}(Optimality: necessary and sufficient conditions.)
 \label{th:linear}
If the function  $\mathcal{H}(t, u, x, p)$ is linear in $(u, x)$  on $\tilde{U} \times X$, then  the conditions (\ref{eq:dHsup}) and (\ref{eq:dHinf})  in  Theorem \ref{th:S} are also necessary.
\end{theorem}

\section{Proofs}
{\sl Proof of Theorem \ref{th:S}}.
%
 Put ${\bf y}:={\bf x}-\hat {\bf x}$. Since 
$$
dy_t=\left(f(t,u_t,x_t)-f(t,\hat x_t,\hat u_t)\right)dt +\left(\sigma(t,u_t,x_t)-\sigma(t,\hat u_t,\hat x_t)\right)dW_t,
$$
it follows by the Ito product formula
\bean
    d( p_t^Ty_t)
    &=& p_t^T\, dy_t +y_t\,d{p}_t^T+h^T_t\big(\sigma(t,u_t,x_t)-\sigma(t,\hat u_t,\hat x_t)\big)dt\\ 
    &=& p_t^T\big(f(t,u_t,x_t)-f(t,\hat u_t,\hat x_t)\big)\,dt + h^T_t(\sigma(t,u_t,x_t)-\sigma(t,\hat u_t,\hat x_t)\big)dt\\ 
    &&   -\cH_x(t,\hat u_t, \hat x_t, p_t^T,h^T_t)y_t\, dt
           + h^T_ty_t\, dW_t\\ && +p_t^T\big(\sigma(t,u_t,x_t)-\sigma(t,\hat u_t,\hat x_t)\big)\,dW_t. 
\eean
Using its integral form we have:  
\bean
   p_T^Ty_T&=&p_0^Ty_0
    + \int_0^T p_t^T\, dy_t +y_t\,d{p}_t^T+\int_0^Th^T_t\big(\sigma(t,u_t,x_t)-\sigma(t,\hat u_t,\hat x_t)\big)\,dt\\ 
    &=& \int_0^T p_t^T\left[\big(f(t,u_t,x_t)-f(t,\hat u_t,\hat x_t)\big) + h^T_t(\sigma(t,u_t,x_t)-\sigma(t,\hat u_t,\hat x_t)\big)\right]dt\\  
    &&   -\int_0^T\cH_x(t,\hat u_t, \hat x_t,  p_t^T,h^T_t)\, y_t\, dt
     + \int_0^T h^T_ty_t\, dW_t\\ && +\int_0^T p_t^T\left(\sigma(t,u_t,x_t)-\sigma(t,\hat u_t,\hat x_t)\right) dW_t. 
\eean
The stochastic integrals 	above have zero expectations due to our
assumptions. Since 
 $p^T_T=0$ and $y_0=0$  we get that 
 \begin{align*}
 &\E \int_0^T p_t ^T\left[\big(f(t,u_t,x_t)-f(t,\hat u_t,\hat x_t)\big) + h^T_t(\sigma(t,u_t,x_t)-\sigma(t,\hat u_t,\hat x_t)\big)\right] dt\\&=\E    \int_0^T\cH_x(t,\hat u_t, \hat x_t,  p_t^T,h^T_t)\, y_t\, dt.
 \end{align*}

With this we  complete the proof easily. Indeed, 

\begin{align*}
&\Delta J_T :=  J_T(\hat{\boldsymbol{u}}) - J_T(\boldsymbol{u})=\E \int\limits_{0}^{T} \big(g( t,\hat u_t, \hat x_t) - g(t,u_t, x_t) \big)\, 
dt\\
&=\E\int_{0}^{T} \big(\cH(t,\hat u_t, \hat x_t, p_t^T,h^T_t) - \cH(t,u_t, x_t, p_t^T,h^T_t)\big)\,dt\\
&\quad - \E \int_0^T p_t^T\left[\big(f(t,u_t,x_t)-f(t,\hat x_t,\hat u_t)\big) + h^T_t(\sigma(t,u_t,x_t)-\sigma(t,\hat u_t,\hat x_t)\big)\right]dt\\ 
&= \E\int_{0}^{T} \big(\cH(t,\hat u_t, \hat x_t,  p_t^T,h^T_t) - \cH(t,u_t, x_t,  p_t^T,h^T_t)\big)\,dt\\
&\quad -\E    \int_0^T\cH_x(t,\hat u_t, \hat x_t,  p_t^T,h^T_t)\,y_t\, dt.
\end{align*}


In virtue of Assumption \ref{a:concave} on concavity 
\begin{align}
&\cH(t,u_t, x_t,  p_t^T,h^T_t) - \cH(t,\hat u_t, \hat x_t, p_t^T,h^T_t) \nonumber\\
    &\le \cH_u( t,\hat u_t, \hat x_t, p_t^T,h^T_t) \big( u_t - \hat u_t \big)
     + \cH_x(t,\hat u_t, \hat x_t,  p_t^T,h^T_t) \big( x_t - \hat x_t \big). \label{conc}
\end{align}
Using  this inequality we obtain  from the above formula the bound 
\beq
\label{bound}
 \Delta J_T \ge \E\int_0^T\cH_u(t,\hat u_t, \hat x_t,  p_t^T,h^T_t)\big(\hat u_t-u_t\big)\, dt. 
\eeq
%

Taking limits we get the inequalities 
 $$
    \limsup\limits_{T \to \infty}\Delta J_T\ge\limsup\limits_{T \to \infty}\E\int_0^T\cH_u(t,\hat u_t, \hat x_t, p_t^T,h^T_t)\big(\hat u_t-u_t\big)\, dt
$$
and
$$
       \liminf\limits_{T \to \infty}\Delta J_T\ge\liminf\limits_{T \to \infty} \E\int_0^T\cH_u(t,\hat u_t, \hat x_t,  p_t^T,h^T_t)\big(\hat u_t-u_t\big)\, dt. 
       $$
 leading to the sufficient conditions of optimality (\ref{eq:dHsup}) and (\ref{eq:dHinf}).
\fdem

\smallskip
{\sl Proof of Theorem \ref{th:linear}}.
Due to assumed linearity the inequalities  (\ref{conc})  and  (\ref{bound}) hold as  equalities and the above arguments   leads to the identity  
\begin{align*}
  \Delta J_T & =\E\int_0^T\cH_u(t,\hat u_t, \hat x_t,  p_t^T,h^T_t)\big(\hat u_t-u_t\big)\, dt.
  \end{align*}
Taking here upper and lower limits as  $T\to\infty$, we get that  (\ref{eq:dHsup}) and (\ref{eq:dHinf})  are necessary conditions for the optimality.  \fdem

\smallskip
\begin{remark}
Note that it may happen that the limit of $ p_t^T$ as $T\to \infty$ is not  equal to $\hat p_t^\infty$ and even does not exist. 
\end{remark}

\section{Examples}
\textbf{Example 1.}
Let ${\bf x}={\bf x}^{\bf u}$ be given by the SDE 
\begin{equation*}
  dx_t= u_t\, dt + \sigma_t dW_t,\qquad x_0=0, 
\end{equation*}
where phase space of controls is $U:=[-1,1]$ and the process $\sigma$ is bounded. Let 
\begin{equation*}
  J_T(\boldsymbol{u}) = \E \int_{0}^{T}  f_t\,x_t\, dt,
\end{equation*}
where $f\ge 0$ is a predictable process in $L^{22}$. 

The Hamilton--Pontryagin function is $ 
\mathcal{H}(t,u, x, p)= p\,u + f_t x + h \sigma_t
$ 
and $\mathcal{H}_x(t,u, x, p)=f_t$.   The BSDE (\ref{eq6}) is reduced to the form 
$$
dp^T_t=-f_tdt+ h^T_tdW_t, \qquad p^T_T=0,  
$$
and admits an ``explicit" expression, modulo the martingale representation theorem. Indeed, according to the latter,  there is a  predictable process ${\bf h}^T\in L^{22}$, uniquely determined, such that 
$$
 \int_0^{T}f_s\,ds=\E \int_0^{T}f_s\,ds +  \int_0^{T}h^T_s\,dW_s. 
$$
It remains to put 
$$
p^T_t:=\E \int_0^{T}f_s\,ds-\int_0^{t}f_s\,ds+\int_0^{t}h^T_s\,dW_s. 
$$ 
Note that 
$$ 
p_t^T = \E\left[ \int_t^{T}f_s\,ds \,\bigg|\, \mathcal{F}_t \right]\ge 0.
$$
 It follows that the control  $\hat {\bf u} \equiv 1$ is overtaking optimal
 because the relation  
\begin{align*}
  &\liminf\limits_{T \to \infty}\E \int_{0}^{T}  p_t^T \left(1 - u_t\right) dt \ge 0
\end{align*}
holds for any admissible  $\boldsymbol{u}$.  By Theorem \ref{th:linear} this condition is necessary and sufficient for the overtaking optimality. 
\medskip

\textbf{Example 2.} 
Let the controlled process   ${\bf k}={\bf k}^{\bf u}$ be   given by the SDE 
$$ 
dk_t = \left(u_t-\delta\,k_t\right)dt +\sigma dW_t, \qquad k_0 > 0,  
$$
$U:=[\underline u,\bar u]$, $\delta>0$. Let 
$$
J_T({\bf u})=\E \int_0^T e^{-rt} \left( \pi_t\,  
k_t - u_t - \frac{1}{2} u^2_t \right) dt , 
$$
where ${\pi}\ge 0$ is a bounded predictable process, $r> 0$.

In the context of corporate finance, ${k_t}$ is a capital with marginal productivity $\pi_t$ and ${u_t}$ is the intensity of  investments at time $t$. Without investments the capital is decreasing with the rate $\delta$. The goal of the investor is to maximize in a long run expected discounted  profits. The quadratic term in the functional is interpreted as the losses for an adaptation of investments.

Puting  $x_t=e^{\delta\,t} k_t$, we solve the SDE and transform the problem to the one with  dynamics depending only on the control: 
$$ 
dx_t = e^{\delta\,t}\,u_t\, dt+\sigma e^{\delta t}dW_t \qquad x_0=k_0 > 0, 
$$
and  the functional 
$$
J_T({\bf u})=\E \int_{0}^{T} e^{-r\,t} \left( e^{-\delta\, t}\,\pi_t\,  
x_t - u_t -  \frac{1}{2} u^2_t \right) dt.
$$

The Hamilton--Pontryagin function has the form 
$$ 
\mathcal{H}(t,u, x, p) = e^{\delta\,t}p\, u + e^{\delta\,t}\sigma h+  e^{-r\,t} \left(e^{-\delta\, t}\, \pi_t\, 
x - u -  \frac{1}{2} u^2\right),
$$
$\mathcal{H}_x(t,u, x, p)= e^{-(r+\delta )t}\pi _t$ and 
$$
  p_t^T := \E \left[ \int_{t}^{T}e^{-\left(r+\delta\right)
s}\pi_s\, ds \,\bigg|\, \mathcal{F}_t \right].
$$
Define the bounded continuous process ${\bf q}$ with 
$$
q_t := e^{(r+\delta)t} \E \left[ \int_{t}^{\infty}e^{-(r+\delta)s
} \, \pi_s\, ds \,\bigg|\, \mathcal{F}_t \right]
$$
(in economics, $q_t$ is called  ``marginal Tobin's $q$"). 

Let us  check that the control  $\hat {\bf u}:={\bf q}-1$ is overtaking optimal. 

For $T\ge t$ we have
using the tower property of expectations that
\bean 
\cH_u (t,\hat u_t,\hat x_t, p_t^T)
& =&
e^{\delta\,t}\hat p_t^T - e^{-r\,t} - e^{-r\,t} \hat u_t
 =  e^{\delta\,t}\big( p_t^T - e^{-(r+\delta)\,t} q_t\big) \\
& =& - e^{\delta\,t}\,\E\big[ e^{-\left(r+\delta\right) T}\,q_T| \cF_t \big]
 = - e^{-r\,T}\,e^{-\delta\left(T-t\right)}\,\E\big [q_T| \cF_t \big].
\eean
Taking into account the boundedness of the control, the condition  (\ref{eq:dHinf}) holds as the equality: 
$$
  \liminf_{T \to \infty} e^{-r\,T}\,\E \int_{0}^{T} e^{-\delta\left(T-t\right)}\,\E\!\left[\left.q_T\right| \mathcal{F}_t \right] \left(u_t-\hat 
 u_t\right) \, dt =0. 
$$

In virtue of Theorem \ref{th:S}   the control $\hat {\bf u}:={\bf q}-1$ is overtaking optimal. 

\section{Conclusions}
We extended a sufficient condition of optimality of the deterministic problem with infinite-time horizon to stochastic optimal control problems. We call it Seierstad condition in order to distinguish it from Mangasaryan and Arrow sufficient conditions which have already their stochastic analouges in the literature, see e.g.~\cite{OksendalSulem}. Seierstad condition become both sufficient and necessary, when Hamilton--Pontryagin function is lenear with respect to control and state variables.






\begin{thebibliography}{plain}

\bibitem{ArkinSaksonov1978} 
Arkin, V.I., Saksonov, M.T.: Necessary optimality conditions in control problems for stochastic differential equations. Dokl. Akad. Nauk SSSR 244(1), 11--15 (1979)
\bibitem{Belyakov2020} 
 Belyakov, A.O.: On sufficient optimality conditions for infinite horizon optimal control problems. Trudy Mat. Inst. Steklova 308,  65--75 (2020)
\bibitem{Bismut1973}
Bismut, J.-M.: Conjugate convex functions in optimal stochastic control. J. Math. Anal. Appl. 44, 384--404 (1973)
\bibitem{Bismut1976a}
Bismut, J.-M.: Linear quadratic optimal stochastic control with random coefficients. SIAM Journal on Control and Optimization 14, 3 (1976)  
\bibitem{Bismut1976}
Bismut, J.-M.: Growth and optimal intertemporal allocation of risks. J. Econ. Theory 10, 239--257 (1975)
\bibitem{CZ} 
Cvitanic, J., Zhang, J.: Contract Theory in Continuous Time Models. Springer, (2012) 
\bibitem{Kabanov1978} 
Kabanov, Yu.M.: On the Pontryagin maximum principle for linear stochastic differential equations. In Probabilistic models and management of economic processes, pp. 65-94. CEMI, Moscow (1978)   
\bibitem{Marinelli}
Marinelli, C.: The stochastic goodwill problem. European Journal of Operational Research 176(1), 389--404 (2007) 
\bibitem{NA}
Nerlove, M.,  Arrow, J. K.: Optimal advertising policy under dynamic conditions. Economica 29, 129--142 (1962)
\bibitem{OksendalSulem}
\O ksendal, B., Sulem, A.: Applied Stochastic Control of Jump Diffusions. Springer International Publishing (2019)
\bibitem{Romer}
 Romer, D.: Advanced Macroeconomics.  McGraw-Hill, California (2012)
\bibitem{PP1990}
Pardoux, E., Peng, S.: Adapted solution of a backward stochastic differential equation. Systems \& Control Letters 28(1),  55--61 (1990) 
\bibitem{P1990}
Peng, S.: A general stochastic maximum principle for optimal control problems. SIAM Journal on Control and Optimization 28(4), 966--979 (1990) 
\bibitem{Seierstad1977} Seierstad, A.: A sufficient condition for control problems with infinite horizons, 
Memorandum from Institute of Economics, University of Oslo (January 12, 1977)
\bibitem{Seierstad1978} 
Seierstad, A., Syds\ae ter, K.:  Optimal Control Theory with Economic Applications. North-Holland, Amsterdam (1987)
\end{thebibliography}
\end{document}